\documentclass{scrartcl}
\usepackage{cite}  
\usepackage{graphicx}  
\usepackage{amsmath}  
\usepackage{amssymb}
\usepackage{url}       
\usepackage{algpseudocode}
\usepackage{enumerate}
\usepackage{units}
\usepackage{amsthm}
\usepackage[english]{babel}

\begin{document}

\title{\LARGE \bfseries A Mixed-Integer QCQP Formulation of AC Power System 
  Upgrade Planning Problems}
\subtitle{Technical Note}
\author{ 
  Sandro Merkli\thanks{
    Automatic Control Lab,
    ETH Zurich, Physikstrasse 3, 8092 Zurich
    \texttt{ merkli@control.ee.ethz.ch}}
}

\def \todo {\textbf{todo:} }
\def \div { \mathrm{div} }
\def \rot { \mathrm{\textbf{rot}} }
\def \minfty { -\infty }
\def \infint { \int_{-\infty}^\infty }
\def \infsum { \sum_{k = -\infty}^\infty }
\def \adj { \mathrm{adj} }
\def \lps { \quad \laplace \quad }
\def \B { \mathbb B }
\def \C { \mathbb C }
\def \N { \mathbb N }
\def \R { \mathbb R }
\def \Z { \mathbb Z }
\def \Q { \mathbb Q }
\def \Ac { \mathcal A }
\def \Bc { \mathcal B }
\def \Dc { \mathcal D }
\def \Ec { \mathcal E }
\def \Fc { \mathcal F }
\def \Gc { \mathcal G }
\def \Xc { \mathcal X }
\def \Ic { \mathcal I }
\def \Zc { \mathcal Z }
\def \Qc { \mathcal Q }
\def \Rc { \mathcal R }
\def \Uc { \mathcal U }
\def \Nc { \mathcal N }
\def \ei { \varepsilon_{\text{int}} }
\def \xtt { x_{t|t}   }
\def \gdw { \;\; \Longleftrightarrow \;\; }
\def \diag { \operatorname{diag} }
\def \minim { \operatorname*{minimize} }
\def \maxim { \operatorname*{maximize} }
\def \st { \operatorname*{subject\ to} }
\def \dom { \operatorname{dom} }
\def \img { \operatorname{im} }
\def \conv { \operatorname{conv} }
\def \cone { \operatorname{cone} }
\def \projop { \operatorname{proj} }
\def \jw  { j\omega }
\def \bmb { \begin{bmatrix} }
\def \bme { \end{bmatrix} }
\newcommand{\cve}[1] {#1}
\newcommand{\jcve}[1] {\bar{#1}}
\newcommand{\drop}[1] {{ }}

\maketitle
\pagestyle{empty}
\thispagestyle{empty}

\section{Introduction}
This document describes the detailed reformulation of a power system upgrade
planning problem into a more generic quadratically constrained quadratic 
problem (QCQP). The problem is one of deciding what lines to upgrade in an
existing power system in order to enable it to handle specific load situations.

The reformulation presented here is based on other formulations already used in
the field, such as the ones presented in~\cite{Jabr2013} and~\cite{Low2014b}. This
document was created as a reference for publications by the author that refer
to the planning problem but avoid outlining reformulation details.

\section{Original formulation}
The original problem is given as follows:
\begin{subequations}
  \label{eqn:upgradeproblem}
  \begin{align}
  \minim_{a,\tilde v^k,\tilde s^k} &\;\; f(a) \\
    \st &\;\; Aa \le b, \\
        &\;\; a \in \{0,1\}^{n_u}, \label{eqn:up_upg}\\ 
        &\;\; Y_{\text{upg}} = Y + \sum_{i=1}^{n_u} (a_i \cdot \delta Y_i), \label{eqn:up_upg2} \\
        &\;\; I_{\max,\text{upg},jl} = I_{\max,jl} + \sum_{i \in \mathcal U_{jl}} 
          (a_i \cdot \delta I_{jl,i}), \label{eqn:up_upg2} \\
        &\;\; \diag(\tilde v^k)\overline{Y_{\text{upg}}\tilde v^k} = \tilde s^k, \label{eqn:up_kir}\\
        &\;\; v_{\min,j} \le |\tilde v^k_j| \le v_{\max,j}, \\
        &\;\; |Y_{\text{upg},jl}||\tilde v^k_j - \tilde v^k_l | \le I_{\max,\text{upg},jl}, \\
        &\;\; s_{\min,j}^k \le \tilde s^k_j \le s_{\max,j}^k, \label{eqn:up_powbounds}\\
        &\;\; \forall k \in \{1,\ldots, K\} \nonumber
  \end{align}
\end{subequations}
where $a \in \{0,1\}^{n_u},Y_\text{upg} \in \C^{N \times N}, \tilde v^k \in
\C^N, \tilde s^k \in \C^N$ are variables and the rest is data. The index $k$
refers to the different scenarios that have to be accounted for. The vector $a$
contains the upgrade decisions, $Y_\text{upg}$ is the upgraded Laplacian matrix
of the power grid, $\tilde s^k,\tilde v^k$ are the new vectors of powers and
voltages, respectively. 

\section{Reformulation}
In the following, the index $k$ will be dropped and only one scenario will be
considered. This is done because it simplifies notation and the reformulation
procedure is the same for each individual scenario. The reformulated problem
will be given as 
\begin{subequations}
  \label{eqn:upgradeproblem2}
  \begin{align}
  \minim_{a,z,y} &\;\; f(a) \\
         \st &\;\; Aa \le b \label{eqn:up2_poly}\\
             &\;\; a \in \{0,1\}^{n_u}\label{eqn:up2_bin} \\ 
             &\;\; \alpha_i \le z^TQ_iz + q_i^Ty + m_i^Ta \le \beta_i\label{eqn:up2_quad} \\
             &\;\; \forall i \in \{1,\ldots, I\},  \nonumber
  \end{align}
\end{subequations}
where the following sections outline variable correspondences and how the
constraints are brought into the standard form above. Note that not all
constraints~\eqref{eqn:up2_quad} have all terms present. For example, some
constraints have $\alpha_i = -\infty$ or $\beta_i = \infty$, meaning they are
only one-sided constraints. This changes nothing about the convexity properties
of the problem (the voltage constraints alone make it non-convex).

\subsection{Variables}
\label{ssec:variables}
The variables of the reformulated problem are as follows:
\begin{itemize}
  \item[-] Variables $v_r \in \R^N,v_q \in \R^N$ representing the real and
    imaginary parts of $\tilde v$,
  \item[-] Variables $l_r \in \R^{2L}$, $l_q \in \R^{2L}$ representing real and
    reactive powers flowing into buses from lines. There are a total of four
    variables per line: Real and reactive powers flowing into the line from
    the two buses the latter connects. The distinction between the two 
    power flow directions is required since there are power losses along the
    line, so power flowing into the line from one bus does not usually 
    equal the power flowing out on the other side.
  \item[-] Variables $a \in \{0,1\}^{n_u}$ are the same as in the original
    formulation.
\end{itemize}
We will use the shorter notations $z := \bmb v_r^T & v_q^T \bme^T$ and $y :=
\bmb l_r^T & l_q^T \bme^T$. If there are multiple snapshots, copies of the
entire $z$ and $y$ vectors along with all the constraint that follow in the
next sections have to be added for each snapshot.

\subsection{Voltage constraints}
A voltage constraint 
\begin{equation} 
  v_{\min,j} \le |\tilde v_j| \le v_{\max,j} 
\end{equation}
can be rewritten in the new variables as 
\begin{equation}
  v_{\min,j}^2 \le v_{r,j}^2 + v_{q,j}^2 \le v_{\max,j}^2.
\end{equation}
and hence as
\begin{equation}
  \label{eqn:voltfin}
  v_{\min,j}^2 \le z^T Q_j z \le v_{\max,j}^2.
\end{equation}
where $Q_j$ has entries $1$ in positions $(j,j)$ and $(N+j,N+j)$ and zeros
everywhere else. Equation~\eqref{eqn:voltfin} is now in the form
of~\eqref{eqn:up2_quad}.

\subsection{Current constraints}
A current constraint as follows:
\begin{equation}
  |Y_{\text{upg},jl}||\tilde v^k_j - \tilde v^k_l | \le I_{\max,\text{upg},jl}, \\
\end{equation}
requires a bit more work to reformulate. We first move everything related
to upgrades to the right-hand side and square both sides:
\begin{equation}
  |\tilde v_j - \tilde v_l |^2 \le \frac{I_{\max,\text{upg},jl}^2}
    {|Y_{\text{upg},jl}|^2}.
\end{equation}
We then rewrite the left-hand side as a function of $x$:
\begin{equation}
  \begin{aligned} |\tilde v_j - \tilde v_l|^2 &= \Re(\tilde v_j-\tilde v_l)^2 
      + \Im(\tilde v_j-\tilde v_l)^2 \\
    &= (v_{r,j} - v_{r,l})^2 + (v_{q,j} - v_{q,l})^2\\
    &= v_{r,j}^2  + v_{r,l}^2 + v_{q,j}^2 + v_{q,l}^2 
      -2 v_{r,j}v_{r,l}-2 v_{q,j}v_{q,l} \\
    &= z^T Q_{jl} z
  \end{aligned}
\end{equation}
where now $Q_{jl}$ is has:
\begin{itemize}
  \item[-] entry $1$ at $(j,j), (l,l), (N+j,N+j), (N+l,N+l)$,
  \item[-] entry $-1$ at $(j,l),(l,j), (N+j,N+l), (N+l,N+j)$,
  \item[-] zeros everywhere else.
\end{itemize}
The right-hand side depends on the upgrade choices. More specifically, we can
write 
\begin{equation}
  Y_{\mathrm{upg},jl} = Y_{jl} + \sum_{i=1}^{n_u}( a_i \cdot (\delta Y_i)_{jl}).
\end{equation}
as well as
\begin{equation}
  I_{\max,\text{upg},jl} = I_{\max,jl} + \sum_{i=1}^{n_u} (a_i \cdot \delta I_{jl,i})
\end{equation}
Since only one upgrade choice can be made for each line, the fraction can be 
rewritten as an equation that is linear in $a$:
\begin{equation}
  \frac{I_{\max,\text{upg},jl}^2}{|Y_{\text{upg},jl}|^2} = 
  \frac{I_{\max,jl}^2}{|Y_{jl}|^2} + \sum_{i \in \mathcal U_{jl}}
    a_i \left[\left(\frac{\delta I_{jl,i} + I_{\max,jl} }
    {\delta (Y_i)_{jl} + |Y_{jl}| } 
    \right)^2 - \frac{I_{\max,jl}^2}{|Y_{jl}|^2} \right].
\end{equation}
Overall, we can now write the line current constraints as 
\begin{equation}
  \label{eqn:currfin}
 x^TQ_{jl}x + m_{jl}^Ta \le u_{jl},
\end{equation}
where $m_{jl}$ and
$u_{jl}$ are defined as
\begin{equation}
  u_{jl} = \frac{I_{\max,jl}^2}{|Y_{jl}|^2},
\end{equation}
and 
\begin{equation}
  (m_{jl})_i = -\left[\left(\frac{\delta I_{jl,i} + I_{\max,jl} }
    {\delta (Y_i)_{jl} + |Y_{jl}| } 
    \right)^2 - \frac{I_{\max,jl}^2}{|Y_{jl}|^2} \right].
\end{equation}
Equation~\eqref{eqn:currfin} now has the same structure as~\eqref{eqn:up2_quad}.

\subsection{Line power constraints}
The concept used to implement the line constraints is similar to the one
in~\cite{Jabr2013}: We implement Big-M type constraints that enforce one out of
several equalities depending on which upgrade is selected.  We start from the
power flow equations:
\begin{equation}
  \begin{aligned}
    &Y_{\text{upg}}(a) = Y + \sum_{i=1}^{n_u} a_i \delta Y_i, \\
    &\diag(\tilde v)\overline{Y_{\text{upg}}(a)\tilde v} = \tilde s.
  \end{aligned}
\end{equation}
To avoid having higher-order products, this can be rewritten into a case
distinction between different cases of the upgrade vector (case 1 being $\hat
a^1$, case 2 being $\hat a^2$ and so on):
\begin{equation}
  \tilde s_j = 
  \begin{cases}
    \tilde v_j \overline{Y_{\text{upg}}(\hat a^1)\tilde v} & \text{if } a = \hat a^1 \\
    \tilde v_j \overline{Y_{\text{upg}}(\hat a^2)\tilde v} & \text{if } a = \hat a^2 \\
    \vdots \\
    \tilde v_j \overline{Y_{\text{upg}}(\hat a^2)\tilde v} & \text{if } a = \hat a^{2^{n_u}}.
  \end{cases}
\end{equation}
While this would work in theory, there would be a number of constraints
exponential in the number of upgrade options. Since not all upgrades affect
all lines, many of these constraints could be omitted. However, the number of
constraints required per bus $j$ would still be exponential in the number of
upgrade options affecting any of the lines connected to bus $j$.\\
In order to further reduce the number of constraints required, we can introduce
the case distinction per line $(j,l)$ (direction matters due to line losses)
instead of per bus $j$:
\begin{equation}
  \tilde s_{jl} = 
  \begin{cases}
    \tilde v_j \overline{Y_{\text{upg},jl}(\hat a^1)(\tilde v_l-\tilde v_j)}, 
      & \text{if } a = \hat a^1 \\
    \tilde v_j \overline{Y_{\text{upg},jl}(\hat a^2)(\tilde v_l-\tilde v_j)}, 
      & \text{if } a = \hat a^2 \\
    \vdots \\
    \tilde v_j \overline{Y_{\text{upg},jl}(\hat a^2)(\tilde v_l-\tilde v_j)},
      & \text{if } a = \hat a^{n_u(jl)}
  \end{cases}
\end{equation}
where $\tilde s_{jl}$ is to be read as ``power flowing into bus $j$, out of the
line $(j,l)$'' and $n_u(jl)$ is the number of upgrades affecting the line
$(j,l)$. The difference to the previous case is that we only have to consider
the upgrade options for that particular line. Additionally, we know only one
upgrade option can be chosen, so we only need a number of constraints linear in
the number of upgrade options per line. The downside is that variables $\tilde
s_{jl}$ have to be introduced as well as additional constraints building the
power balance for each bus:
\begin{equation}
  \label{eqn:powerbalance}
  s_{\min,j} \le \sum_{l} \tilde s_{jl} \le s_{\max,j}, 
\end{equation}
where $\tilde s_{jl}$ is not added as an additional variable but written with
the help of the $l_r$ and $l_q$ in $y$ (see Section~\ref{ssec:variables} for an
overview of the variables):
\begin{equation}
  \tilde s_{jl} = (l_r)_{jl} + \sqrt{-1} (l_q)_{jl}.
\end{equation}
In the above, the index $jl$ is used to refer to the entries in $l_r$ and $l_q$
that are assigned to the real and reactive parts of the power flowing into bus
$j$ from the line between buses $j$ and $l$.  As for the actual implementation,
for a line $(j,l)$ and upgrade option $k$ affecting it, we would have the
constraints
\begin{equation}
  \label{eqn:linepow_lin1}
  \begin{aligned}
    &-M(1-a_k) \le (l_r)_{jl} - \Re\left(
    \tilde v_j \overline{Y_{\text{upg},jl}(a_k = 1)(\tilde v_l-\tilde v_j)}
    \right) \le M(1-a_k), \\
    &-M(1-a_k) \le (l_q)_{jl} - \Im\left(
    \tilde v_j \overline{Y_{\text{upg},jl}(a_k = 1)(\tilde v_l-\tilde v_j)}
    \right) \le M(1-a_k),
  \end{aligned}
\end{equation}
where $M$ is large enough that if $a_k = 0$, the two constraints can never be
active. In addition to the above, a ``no upgrade case'' has to be added as
well: The case for which no upgrade was performed to this line. This is
implemented by the constraints
\begin{equation}
  \label{eqn:linepow_lin2}
  \begin{aligned}
    &-M\sum_k a_k \le (l_r)_{jl} - \Re\left(
    \tilde v_j \overline{Y_{jl}(\tilde v_l-\tilde v_j)}
    \right) \le M\sum_k a_k, \\
    &-M\sum_k a_k \le (l_q)_{jl} - \Im\left(
    \tilde v_j \overline{Y_{jl}(\tilde v_l-\tilde v_j)}
    \right) \le M\sum_k a_k,
  \end{aligned}
\end{equation}
where the summation over $k$ only includes the $k$ that affect the given line.
Writing the real and reactive parts as functions of $v_r, v_q$ leads to:
\begin{equation}
  \label{eqn:linepow_quad}
  \begin{aligned}
  \Re\left( \tilde v_j \overline{Y_{jl}(\tilde v_l-\tilde v_j)} \right)
    &= -\Re(Y_{jl})(v_r)_j^2 + \Re(Y_{jl})(v_r)_j(v_r)_l - \Im(Y_{jl}) (v_r)_j(v_q)_l \\
    &\quad + \Re(Y_{jl})(v_q)_j(v_q)_l - \Re(Y_{jl}) (v_q)_j^2 + \Im(Y_{jl})(v_q)_j(v_r)_l, \\
  \Im\left( \tilde v_j \overline{Y_{jl}(\tilde v_l-\tilde v_j)} \right) 
    &= \Re(Y_{jl})(v_q)_j(v_r)_l - \Im(Y_{jl})(v_q)_j(v_q)_l + \Im(Y_{jl})(v_q)_j^2 \\
    &\quad - \Re(Y_{jl})(v_r)_j(v_q)_l - \Im(Y_{jl})(v_r)_j(v_r)_l + \Im(Y_{jl})(v_r)_j^2.
  \end{aligned} 
\end{equation}
A similar procedure is used for~\eqref{eqn:linepow_lin1}, simply by replacing
$Y_{jl}$ with $Y_{\mathrm{upg},jl}(a_k = 1)$.
This means the line power constraints are now quadratic in $z$ and linear in
$y$ and $a$, as was required. The data for the $Q$ terms is given
in~\eqref{eqn:linepow_quad}. The data for the $q$ terms is determined by how
the variables are sorted in $l_r$ and $l_q$ and
by~\eqref{eqn:linepow_lin1}-\eqref{eqn:linepow_lin2}.  The data for the $m$
terms is also determined by~\eqref{eqn:linepow_lin1}-\eqref{eqn:linepow_lin2}.
In summary, the constraints in~\eqref{eqn:powerbalance} combined with those
in~\eqref{eqn:linepow_lin1}-\eqref{eqn:linepow_lin2} implement the
constraints~\eqref{eqn:up_upg2},~\eqref{eqn:up_kir}
and~\eqref{eqn:up_powbounds} from the original formulation.
Constraints~\eqref{eqn:powerbalance} are easily brought into the form
of~\eqref{eqn:up2_quad} because they are linear constraints. One set of
constraints in~\eqref{eqn:linepow_lin1} (the same holds
for~\eqref{eqn:linepow_lin2}) translates into 4 constraints of the
form~\eqref{eqn:up2_quad} because the different terms on the upper and lower
bounds in~\eqref{eqn:linepow_lin1}-\eqref{eqn:linepow_lin2} lead to two
different $m$ terms.

\bibliographystyle{ieeetr}
\bibliography{biblio}

\end{document}